\documentclass{qtds}

\usepackage[T1]{fontenc}
\usepackage[english]{babel}
\usepackage{amssymb,amsmath,epsfig,graphics,psfrag}
\usepackage[all]{xy}
\begin{document}

%






\def\RR{{\mathbf{R}}}
\def\CC{{\mathbf{C}}}
\def\HH{{\mathbf{H}}}
\def\NN{{\mathbf{N}}}
\def\ZZ{{\mathbf{Z}}}
\def\cH{{\mathcal{H}}}
\def\cR{{\mathcal{R}}}
\def\cP{{\mathcal{P}}}
\def\cD{{\mathcal{D}}}
\def\cI{{\mathcal{I}}}
\def\cS{{\mathcal{S}}}
\def\cM{{\mathcal{M}}}

\authorrunninghead{Emmanuel Dufraine}
\titlerunninghead{homotopy classes of vector fields}

\title{About homotopy classes of non-singular \\ vector fields on the three-sphere}


\author{Emmanuel Dufraine\thanks{This work was partially supported by E.S.F. PRODYN}}
\affil{Mathematics Institute\\ 
University of Warwick\\ 
Coventry CV4 7AL -- U.K.}

\email{dufraine\@@maths.warwick.ac.uk}
\vspace{24pt}
{\begin{minipage}{24pc}
\footnotesize{\ \ \ \ Generically, the set of points along which two non-singular vector fields on the three-sphere are positively (resp. negatively) collinear form a link. We prove that the two vector fields are homotopic if and only if the linking number of those links is zero. We use this criterion to give a new proof of a result of Yano: every non-singular vector field on the three-sphere is homotopic to a non-singular Morse-Smale vector field.}
\end{minipage}}
\vspace{10pt}


\begin{article}

\section*{Introduction}

The study of surface homeomorphisms up to isotopy, following Nielsen and Thurston, has led to many exciting results. For example, the
 work of Franks on area-preserving diffeomorphisms of the annulus, of Handel about fixed point of planar homeomorphisms, of Bestvina \& Handel for surface homeomorphisms, of Gambaudo, van Strien \& Tresser and Llibre \& MacKay about the forcing problem \dots The cornerstone of this theory is the classification result of Thurston (see~\cite{th88}). In~\cite{ma}, MacKay proposes to study the same problem for non-singular vector fields on three-dimensional manifolds. The main goal of this project is to understand the relationship between the geometry of the manifold (in the sense of Thurston,~\cite{th82}) and the possible vector fields in each homotopy class (homotopy via non-singular vector fields). The problem of homotopy of non-singular vector fields has already been studied by (among others) Asimov, Yano, Gompf and G. Kuperberg~\cite{as,ya2,go,kus}. In the last two papers, Gompf and G. Kuperberg give a complete classification of non-singular vector fields on a three-manifold in term of Euler class, spin structure, \dots Those studies were done using a trivialization of the tangent bundle of the manifold; it seems that they lead to very poor results on the dynamical point of view.

One key point of the programme of Thurston is to find a "simplest representative" in each isotopy class. For three-dimensional vector fields, as remarked by MacKay, the notion of "simplest" is still unclear for one could not expect to define this notion with respect to the topological equivalence of vector fields. Indeed, the work of K. Kuperberg on the Seifert conjecture (see \cite{kum,gh}) implies that each non-singular vector field on a three-manifold is homotopic to a vector field without periodic orbit. However, it is still interesting to know if there is an homotopy to a well-known class of vector fields (e.g. Morse-Smale, volume preserving, pseudo-Anosov flows, \dots). One early work in this direction, after the work of Asimov, is the work of Yano, \cite{ya2}, where conditions for the existence of a Morse-Smale vector field in a given homotopy class on a graph manifold (or on a manifold that admits a Morse-Smale vector field) are given in term of Euler class.

The situation on the three-sphere is in many ways simpler than on other three-manifolds. We will study the problem of homotopy on the sphere, expecting to generalise our results and methods to arbitrary three-manifolds.

One purpose of this paper is to give a criterion to decide whether two non-singular vector fields are homotopic or not, directly computable from the vector fields. Generically, the set of points on which two non-singular vector fields $X$ and $Y$ on the three-sphere are positively (resp. negatively) collinear form a link $C_+$ (resp. $C_-$). The linking number of $C_+$ with $C_-$ is well-defined up to its sign.

{\bf Criterion:} {\em the vector fields $X$ and $Y$ are homotopic if and only if the absolute linking number of $C_+$ with $C_-$ is zero.}

We show that the absolute linking number of $C_+$ and $C_-$ is the distance in homotopy classes between $X$ and $Y$.

We expect this criterion to be extended to arbitrary three-manifold, using previous work of Kuperberg and Gompf, \cite{kus,go} and the extension of Hopf's ideas by Pontryagin.

The last part of the paper is devoted to the study of Morse-Smale vector fields on the three-sphere. We show that every vector field on the sphere is homotopic to a Morse-Smale one, using the criterion. For that purpose, it is sufficient to find a non-singular Morse-Smale vector field in each homotopy class. This has already been done by Yano, \cite{ya2}, using action groupe arguments.
Our construction is more geometric than the construction of Yano. The dynamic of the Morse-Smale vector fields is given by a vector field on the two-sphere and is known explicitly.

\section{Definitions and first results}
\subsection{Hopf fibration and Hopf invariant}\label{ss.hopf}
In this section we give the classical definition of the Hopf fibration of the 3-sphere, we give an alternative way to obtain this fibration and finally, we recall the definition of the Hopf invariant.

We consider $S^3$ as the unit sphere of $\RR^4$ with the standard metric: $S^3=\{(x_1,x_2,x_3,x_4)|\ x_1^2+x_2^2+x_3^2+x_4^2=1\}$. Identifying $\RR^4$ with $\CC^2$, by $(z_1,z_2)=(x_1+i x_2,x_3+i x_4)$, the sphere $S^3$ is the set $\{ |z_1|^2+|z_2|^2=1\}$. Each complex line of $\CC^2$ intersects $S^3$ along a great circle. The union of those circles is the 3-sphere and if two of those circles have a non-empty intersection, they are equal. We obtain a fibration from $S^3$ to $\CC P(1)\simeq S^2$ which associate to each point of the 3-sphere the  direction of the corresponding complex line in $\CC P(1)$.

\begin{definition}
The map described above, $\cH \colon S^3 \to S^2$, is called the {\em Hopf fibration}.
\end{definition}

Let us now identify $\RR^4$ with $\HH$, the field of quaternions.
With this identification, the canonical basis of $\RR^4$ is denoted by $\{1,i,j,k\}$. We identify $\RR\cdot 1$ with $\RR$ and $\RR i + \RR j + \RR k$ with $\RR^3$; this subset of $\HH$ is called the set of pure quaternions. One can write every quaternion $q=\cR(q)+\cP(q)$ with $\cR(q)\in\RR$ and $\cP(q)\in \RR^3$, and define the conjugate $\overline q=\cR(q)-\cP(q)$. We obtain the canonical norm on $\RR^4$: $\|q\|=\sqrt{q\overline q}$.

The 3-sphere is the set $S^3=\{q|\ \|q\|=1\}$ and it has a multiplicative group structure, induced by the multiplication of quaternions. For a given point $s$ in $S^3$, one can consider the inner product: $\rho'_s\colon S^3\to S^3$, defined by $\rho'_s(q)=sqs^{-1}$.

\begin{proposition}[\cite{be}, Corollary 8.9.3]\label{p.quaternions}
Let $s$ be a point of $S^3$, the map $\rho'_s$ leaves $\RR^3$ invariant and its restriction $\rho_s=\rho'_s|_{\RR^3}$ belongs to $SO(3)$.

Moreover the map $s\mapsto \rho_s$ from $S^3$ to $SO(3)$ is a surjective group homomorphism with kernel $\{\pm 1\}$. In particular, $SO(3)$ and $\RR P(3)$ are isomorphic.
\end{proposition}

Let $\star$ be a point of $S^2$, the unit sphere of $\RR^3$, we define the map 
\begin{eqnarray*}
 f \colon & S^3 & \to  S^2 \\
   & s  &\mapsto  \rho_s(\star).
\end{eqnarray*}

\begin{lemma}\label{l.hopf}
For $\star=(0,1,0,0)$, $f$ is the Hopf fibration defined above.
\end{lemma}
\begin{proof}
Let $y$ be a point of $S^2$, we have to prove that $f^{-1}(y)$ is a circle included in a complex line of $\CC^2$. Proposition 8.9.4 in~\cite{be} asserts that given $s=\alpha+t$ in $S^3$, with $t \in \RR^3\setminus{0}$ and $\alpha\in \RR$, the axis of the rotation $\rho_s$ is the line $\RR t$, and its angle $\theta$ is given by the relation $\tan \frac \theta 2 = \frac{\|t\|}{|\alpha|}$ if $\alpha \neq 0$, and $\theta  = \pi $ if $\alpha = 0$.

For $u\in[0,2 \pi]$, let us define $c_u=(\cos (u), \sin(u), 0, 0)$ a point in $S^3$. The rotations associated to $c_u$ by the Proposition~\ref{p.quaternions} are the only ones keeping $\star$ invariant. 

Let us denote $y=(0,a,b,c)$, with $a^2+b^2+c^2=1$ and let $s_{y}$ be the point of $S^3$ with coordinates $(a,0,c,-b)$. We have $\rho_{s_{y}}(y)=\star$, therefore $\rho_{c_u}\circ\rho_{s_{y}}(y)=\rho_{c_u\cdot s_{y}}(y)=\star$. Hence,  we get $$f^{-1}(y)=\{c_u\cdot s_{y},u \in [0,2\pi]\}.$$

Let $u\in[0,2\pi]$, $c_u\cdot s_{y}=(a\cos(u),a\sin(u),b\sin(u)+c\cos(u),c\sin(u)-b\cos(u))$, therefore, $f^{-1}(y)$ is circle, invariant by the antipodal map and included in the complex line of slope $\frac{c-i b}a$.
\end{proof}

\begin{remark}
For a different choice of point $\star$, we obtain a new map $f'\colon S^3\to S^2$, and there exists a diffeomorphism $g$ of the 2-sphere such that $f'=g\circ\cH$ (the map $g$ could be the restriction of an isometry of $\RR^3$ to $S^2$).
\end{remark}

Let $f$ be a smooth map from $S^3$ to $S^2$ (we assume that the spheres are oriented) and let $y$ and $z$ be two regular values of $f$. The sets $f^{-1}(y)$ and $f^{-1}(z)$ are oriented links (a link is disjoint union of knots) and the linking number: $\mbox{link}(f^{-1}(y),f^{-1}(z))$ does not depend on the choice of $y$ and $z$.

\begin{definition}
The {\em Hopf invariant} of $f$ is $H(f)=\mbox{\em link}(f^{-1}(y),f^{-1}(z))$.
\end{definition}

The Hopf invariant does not depend on the choice of $f$ in its homotopy class.

\begin{example}
The Hopf invariant of the Hopf fibration is $1$.

Let $f$ be a map from $S^3$ to $S^2$, $h\colon S^3 \to S^3$ and $g\colon S^2 \to S^2$ we have the equalities: 
$$H(f\circ h)=\deg(h)H(f)$$
$$H(g\circ f)=\deg^2(g) H(f).$$
 
We deduce from the second formula that the Hopf invariant of any map $f$ as above (see before lemma~\ref{l.hopf}) is $1$.
\end{example}

\subsection{Combings and framings}
Let us define some geometrical objects on an arbitrary compact oriented 3-manifold $M$. In this section we follow~\cite{wi1,kus,go,gagh}. We denote by $TM$ the tangent bundle of $M$.

\begin{definition}
A {\em combing} of $M$ is a section of the tangent bundle $TM$. 

A {\em framing} of $M$ is a trivialization of the tangent bundle into a product: $TM \simeq M\times \RR^3$.
\end{definition}

A combing is a non-singular vector field on $M$. As every compact oriented 3-manifold has a vanishing Euler characteristic, there always exists a combing on $M$. Every combing is homotopic to a unitary combing for a given Riemannian metric on $M$.

A framing consists of three linearly independent combings whose orientation gives the orientation of the manifold. We remark that with a Riemannian metric on $M$, two linearly independent combings are sufficient to define a framing and every framing is homotopic to an orthonormal framing. It is more difficult to show that there always exists a framing on a three-manifold. We will restrict ourselves to the case $M=S^3$ which is a Lie group, therefore it is easy to see that $S^3$ admits a framing. The next lemma is classical.

\begin{lemma}
On the three-sphere, one can always complete a combing into a framing. Moreover, two such framings are homotopic through framings which complete the initial combing.
\end{lemma}

Therefore, to each non-singular vector field $X$ on $S^3$, we associate a framing $\tau_{_X}$ which is well-defined up to homotopy.

Let $X$ and $Y$ be two unitary non-singular vector fields on the three-sphere and denote by $\tau_{_X}$ and $\tau_{_Y}$ respectively some associated orthonormal framings. We have the following.
\begin{itemize}
\item With respect to $\tau_{_X}$, the vector field $Y$ is a map $Y_{\tau_{_X}}$ from $S^3$ to $S^2$, we can associate to this map its Hopf invariant: $H_X(Y)$;
\item at each point $x$ in $S^3$, we can associate a unique linear map of $SO(3)$,  denoted by $(\tau_{_X}-\tau_{_Y})(x)$, mapping $X(x)$ to $Y(x)$ and the two others combings of $\tau_{_X}$ to the combings of $\tau_{_Y}$. The map $(\tau_{_X}-\tau_{_Y})$ from $S^3$ to $SO(3)$ has a well-defined degree: $[\tau_{_X}-\tau_{_Y}]$.
\end{itemize}

\begin{remark}
As $SO(3)$ is homeomorphic to $\RR P(3)$ and $S^3$ is simply connected, $[\tau_{_X}-\tau_{_Y}]$ is an even integer (see proof of lemma~\ref{l.calcul}).
\end{remark}

\begin{lemma}\label{l.calcul}
The quantity $H_X$ gives an isomorphism between the homotopy classes of non-singular vector fields on the three-sphere and the integers.

Two non-singular vector fields on the three-sphere, $X$ and $Y$, are homotopic if and only if $H_X(Y)=0$.

Moreover, we have the formula: $H_X(Y)=\frac12[\tau_{_X}-\tau_{_Y}]$.
\end{lemma}
\begin{proof}
Let $X$ be a non-singular vector field on $S^3$ and $\tau_{_X}$ an associated framing. The homotopy class of a non-singular vector field $Y$ is uniquely determined by the homotopy class of $Y_{\tau_{_X}}\colon S^3\to S^2$. The Hopf invariant gives an isomorphism between $\pi_3(S^2)$ and $\ZZ$ then it gives an isomorphism between the homotopy classes of non-singular vector fields on $S^3$ and $\ZZ$. 

A continuous map from $S^3$ to $S^2$ is homotopic to a constant if and only if its Hopf invariant is zero. The map $X_{\tau_{_X}}$ is constant, then $X$ and $Y$ are homotopic if and only if $H_X(Y)=0$.

Let us prove now the last formula. Let $\star$ be a point of $S^2$ and define the map $M \colon SO(3)\to S^2$ which associate to a matrix $A$ the point $A(\star)$ of $S^2$. Moreover, denote by $p\colon S^3 \to \RR P(3)\simeq SO(3)$ the standard projection. With $X_{\tau_{_X}}\equiv \star$, there exists $f\colon S^3\to S^3$ such that we have the following commutative diagram.

$$
\xymatrix{&{S^3}\ar[ld]_f\ar[d]|{{\tau_{_X}-\tau_{_Y}}}\ar[r]^{Y}&{S^2}\\
	S^3\ar[r]^p&{SO(3)}\ar[ur]_M&	}
$$

The projection $p$ is of degree two then the integer $[\tau_{_X}-\tau_{_Y}]$ is even and the degree of $f$ is exactly $\frac12[\tau_{_X}-\tau_{_Y}]$. 

Using this diagram, we see that $Y=M\circ p \circ f$ then $H_X(Y)=H(M\circ p \circ f)= \deg(f)H(M\circ p)$. The map $M\circ p$ is exactly the map constructed in section~\ref{ss.hopf} and, $H(M\circ p)=H(\cH)=1$. Therefore we have $H_X(Y) = \deg(f)= \frac12[\tau_{_X}-\tau_{_Y}]$.
\end{proof}

\begin{remark}\label{r.baddef}
The isomorphism we obtain is not well-defined, it depends on the choice of a preferred combing (or framing) on the three-sphere.
\end{remark}

\begin{remark}
If we exchange $X$ and $Y$ we get: $H_X(Y)=-H_Y(X)$.

If $X$, $Y$ and $Z$ are three non-singular vector fields on $S^3$, we have the following equality, for $x$ in $S^3$: $(\tau_{_X}-\tau_{_Y})(x)\cdot (\tau_{_Y}-\tau_{_Z})(x)= (\tau_{_X}-\tau_{_Z})(x)$.

Hence, $[\tau_{_X}-\tau_{_Y}]+[\tau_{_Y}-\tau_{_Z}]=[\tau_{_X}-\tau_{_Z}]$ and $H_X(Y)+H_Y(Z)=H_X(Z)$.
\end{remark}

These remarks lead to the following definition.
\begin{definition}
The {\em distance in homotopy classes} between $X$ and $Y$, denoted by $\cD(X,Y)$, is the absolute value $|H_X(Y)|=|H_Y(X)|$. 
\end{definition}

Following Theorem 5.4 and Theorem A.4 of~\cite{ep}, a diffeomorphism of the three-sphere is either isotopic to the identity (orientation-preserving) or isotopic to \mbox{$R(x_1,x_2,x_3,x_4)=(x_1,x_2,x_3,-x_4)$} (orientation-reversing). We will see in section~\ref{ss.examples} that the vector field $\cH_+$, tangent to the Hopf fibration, satisfies $\cD(\cH_+,R_{\star}(\cH_+))=1$.

\begin{definition}
We associate to a vector field $X$ on the three-sphere its {\em homotopy number}: $$\cI(X)=\frac{(\cD(X,\cH_+)+\cD(X,R_\star(\cH_+))-1)}{2}.$$
\end{definition}

It is clear that $\cD(X,Y)=\cD(g_\star(X),g_\star(Y))$ for any vector fields $X$ and $Y$, and for any diffeomorphism $g$. Therefore we have the formula $\cI(X)=\cI(R_\star(X))$. The next lemma is straightforward. 

\begin{lemma}\label{l.isomorphism}
The homotopy number gives a well-defined isomorphism between the homotopy classes (up to diffeomorphism) of non-singular vector fields on the three-sphere and the natural integers $\NN$.
\end{lemma}

In particular, we get the formula: $\cD(X,R_\star(X))=2\cI(X)+1$ and we obtain that there is a unique homotopy class (up to diffeomorphism) such that $\cD(X,R_\star(X))=1$.

\begin{remark} If we pay attention to the homotopy classes of vector fields up to diffeomorphism, the not well-defined isomorphism of remark~\ref{r.baddef} become the well-defined isomorphism of lemma~\ref{l.isomorphism}.
\end{remark}

\subsection{Intrinsic definition of~$\cD(X,Y)$}

Given two vector fields $X$ and $Y$, it is not easy to compute $H_X(Y)$, because we have to complete $X$ into a framing and to express $Y$ in the coordinates of this framing. We will investigate another way, more direct, to compute $\cD(X,Y)$.

Let $X$ and $Y$ be two non-singular vector fields on the three-sphere. We denote by $C_+$ (resp. $C_-$) the set of points of $S^3$ where $X$ and $Y$ are positively (resp. negatively) collinear.
$$C_+=\{x\in S^3/ X(x)=\lambda Y(x)\mbox{ with } \lambda > 0\}$$
$$C_-=\{x\in S^3/ X(x)=\lambda Y(x)\mbox{ with } \lambda < 0\}$$

\begin{lemma}\label{l.circles}
Generically, $C_+$ and $C_-$ are embedded links in $S^3$ (possibly empty).
\end{lemma}

\begin{proof}
The vector fields $X$ and $Y$ are homotopic to unitary vector fields for a given Riemannian metric. Fixing any orthonormal framing $\tau$ of $S^3$, the vector fields $X$ and $Y$ give a map $(X,Y)$ from $S^3$ to $S^2\times S^2$. The diagonal $D_+$ (resp. the anti-diagonal $D_-$) of $S^2\times S^2$ is the subset of points of the form $(x,x)$ (resp. $(x,-x)$) with $x$ in $S^2$. The set $C_+$ (resp. $C_-$) is the reciprocal image $(X,Y)^{-1}(D_+)$ (resp. $(X,Y)^{-1}(D_-)$).

The spheres $D_+$ and $D_-$ are codimension $2$ submanifolds of $S^2\times S^2$. The transversality theorem asserts that a small perturbation of every map $g\colon S^3\to S^2\times S^2$ is transverse to $D_+$ and $D_-$. Therefore a small perturbation of $X$ and $Y$ makes the map $(X,Y)$ transverse to $D_+$ and $D_-$. Hence $(X,Y)^{-1}(D_+)$ and $(X,Y)^{-1}(D_-)$ are either empty or compact submanifolds of $S^3$, of codimension 2. Then if $C_+$ and $C_-$ are non-empty, they are disjoint unions of embedded circles in $S^3$.
\end{proof}
\begin{remark}
Lemma~\ref{l.circles} is true on arbitrary three-manifold.

Although we made use of a framing in the previous proof, the sets $C_+$ and $C_-$ do not depend on a framing (indeed one can prove lemma~\ref{l.circles} by a local argument).
\end{remark}

Let $M$ be a compact three-manifold, and $f$ and $g$ be two maps from $M$ to $S^2\subset\RR^3$. If the points $f(x)$ and $g(x)$ are never antipodal, that is if the set $\{ x\in M|\ f(x)=-g(x)\}$ is empty or alternatively if $\|f(x)-g(x)\|<2 \pi$ for all $x$ in $S^2$, then $f$ and $g$ are homotopic. This leads to the following lemma.

\begin{lemma}\label{l.vide}
If either $C_+$ or $C_-$ is empty, then $X$ and $Y$ are homotopic.
\end{lemma}

\begin{remark}
The previous argument shows that a perturbation of a vector field leading to a vector field in the same homotopy class can be very large.
\end{remark}

Our goal is to compute the linking number between $C_+$ and $C_-$. For that purpose, we need to give them an orientation; in fact we will give two possible orientations of those links. 

Let us denote by $C_0$ one component of $C_+$, and let $D$ be a small disc, transverse to $X$ and $Y$ at a point $x$ of $C_0$. The vector $X(x)$ gives an orientation of this disc. We consider the Gauss map from $D$ to $S^2$ (with a fixed orientation) which associates to a point of $D$ the direction of the vector field $Y$ at this point.
We oriente $C_0$ such that the orientation of $D$ and the orientation of $C_0$ at $x$ gives the orientation of $S^3$ if the Gauss map is orientation-preserving or gives the opposite of the orientation of $S^3$ otherwise. This choice of orientation does not depend on the choice of $D$ or of $x$.
We oriente this way each component of $C_+$ and $C_-$.

If we take the Gauss map associated to $X$ instead, we obtain the opposite orientation on all the components. We recall that the linking number with the empty set is zero.

\begin{definition}
The {\em absolute linking number} of $C_+$ and $C_-$ is the absolute value of the linking number of $C_+$ and $C_-$ with one of the above orientations: $|\mbox{\em link}(C_+,C_-)|$.
\end{definition}

\begin{remark}
If $C_+$ and $C_-$ are linked by only two components (one for each link), any orientation of those links would give the same linking number, up to sign.

If $C_+$ and $C_-$ are orbits of $X$ (and of $Y$), all their components are oriented as orbits of $X$ (or of $Y$).
\end{remark}

The intrinsic definition of $\cD(X,Y)$ is given by the following.

\begin{lemma}\label{l.critere}
The distance in homotopy classes between two vector fields $X$ and $Y$ is given by the absolute linking number of $C_+$ and $C_-$.
 $$|\mbox{\em link}(C_+,C_-)|=\cD(X,Y)$$
\end{lemma}

\begin{proof}
Let $\tau_{_X}$ be a framing associated to $X$, and $Y_{\tau_{_X}}$ be the vector field $Y$ expressed in this framing. Let $\star$ be the point of $S^2$ such that $X_{\tau_{_X}}\equiv \star$. We can perturb $Y$, staying in the same homotopy class, such that $\star$ and $-\star$ are regular values of $Y_{\tau_{_X}}$.

As $C_+$ is $Y^{-1}_{\tau_{_X}}(\star)$ and $C_-$ is $Y^{-1}_{\tau_{_X}}(-\star)$, we can orientate $C_+$ and $C_-$ with respect to the orientation given by the map $Y^{-1}_{\tau_{_X}}$ and the orientation of $S^3$ and $S^2$. This orientation fits with one of the orientations of $C_+$ and $C_-$ given above, depending on the orientation of $S^2$. By the definition of the Hopf invariant, we have $\mbox{ link}(C_+,C_-)=H_X(Y)$. 

The choice for the orientation of $C_+$ and $C_-$ is not canonical, as well as the choice of orientation of $S^2$, we obtain the result: $|\mbox{link}(C_+,C_-)|=\cD(X,Y)$.
\end{proof}

We obtain the criterion stated in the introduction as a corollary of the previous lemma and lemma~\ref{l.calcul}.

\subsection{First examples}\label{ss.examples}
In this section, we use our criterion to show that $\cH_+$ and $R_\star(\cH_+)$ are in adjacent homotopy classes, we give examples of vector fields homotopic to $\cH_+$ and $R_\star(\cH_+)$ and finally we give a construction of a vector field in each homotopy class.

Let $\cH\colon S^3 \to S^2$ be the Hopf fibration. Choosing an orientation on $S^3$ and on $S^2$ gives a natural orientation of the fibres. We obtain this way a vector field $\cH_+$, tangent to the fibres of $\cH$, such that the linking number between two different orbits is $+1$. The equations of this vector field are: $\cH_+(z_1,z_2)=(iz_1,iz_2)$ if we identifie $\RR^4$ with $\CC^2$ as in the first section. 

In a similar way, we define $\cH_-$ such that the linking number of two orbits is~$-1$: $\cH_-(z_1,z_2)=(iz_1,i\overline{z_2})$. We remark that $\cH_-=R_\star(\cH_+)$. The following lemma was proved with an other method in~\cite{wi1}.

\begin{figure}[htb]
\psfrag{H+}{$\cH_+$}\psfrag{H-}{$\cH_-$}
\centerline{\includegraphics{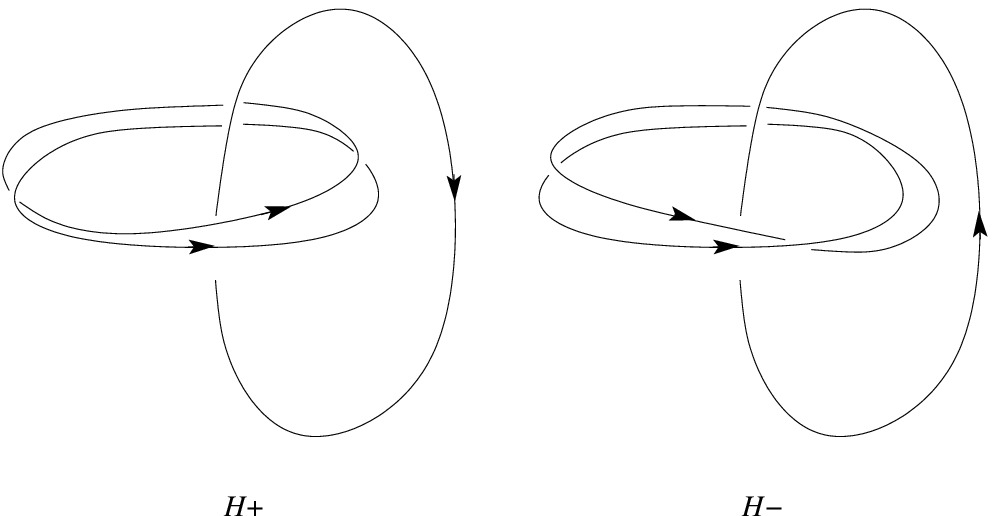}}
\caption{$\cH_+$ and $\cH_-$\label{f.hplusmoins}}
\end{figure}
\begin{lemma}
The two vector fields above, $\cH_+$ and $\cH_-$, lie in adjacent homotopy classes: \mbox{$\cD(\cH_+,\cH_-)=1$}.
\end{lemma}

\begin{proof}
The vector fields $\cH_+$ and $\cH_-$ are positively tangents along $C_+$, the unit circle of $\CC\times\{0\}$ in $\CC^2$ and they are negatively tangent along $C_-$, the unit circle of $\{0\}\times\CC$. The circles $C_+$ and $C_-$ form a Hopf-link, their absolute linking number is 1. Using lemma~\ref{l.critere}, we have $\cD(\cH_+,\cH_-)=1$.
\end{proof}

As a corollary, we obtain that $\cI(\cH_+)=\cI(\cH_-)=0$.

\begin{definition}
Given two relatively prime numbers $p\neq 0$ and $q>0$, a $(p,q)-$Seifert fibration of the three-sphere is a map from $S^3\subset\CC^2$ to $S^2\simeq\CC P(1)$ defined by:
$$\cS_{p,q}\colon (z_1,z_2)\mapsto \left[\frac{z_2^q}{z_1^p}\right]$$
\end{definition}

A Seifert fibration is a continuous map from the three-sphere to the two-sphere, such that the reciprocal image of each point is a circle. Each of those circles are $(p,q)-$knots of the tori $T_{a,b}=\{(z_1,z_2)|\ |z_1|=a,\ |z_2|=b\}$ except the reciprocal image of the North and the South poles which are the cores of those tori and give a Hopf link. Let us denote by $H_{p,q}$ the unitary vector field tangent to the fibres of the $(p,q)-$Seifert fibration. We remark that $H_{1,1}=\cH_+$ and $H_{-1,1}=\cH_-$.

\begin{lemma}[\cite{wi1}, corollary 2.4]\label{l.plusmoins}
The vector fields $H_{p,q}$ are homotopic to $\cH_+$ if $p$ is positive, and to $\cH_-$ if $p$ is negative.
\end{lemma}
\begin{proof}
This is an easy corollary of lemma~\ref{l.vide}.
For $p>0$, the vector fields $H_{p,q}$ and $\cH_+$ are never negatively tangents and it is the same for $H_{p,q}$ and $\cH_-$ whenever $p<0$.
\end{proof}

Let us consider the link consisting of two regular fibres $L_1$ and $L_2$ of a $(p,q)-$Seifert fibration. The absolute linking number between $L_1$ and $L_2$ is the product $|pq|$. We will construct a vector field $X_n$ for each $n$ in $\NN\setminus\{0\}$ which is positively tangent to $H_{n,1}$ along $L_1$ and negatively tangent to $H_{n,1}$ along $L_2$ (and the vector fields are transverse everywhere else). Such a vector field satisfies: $\cD(X_n,H_{n,1})=\cD(X_n,\cH_+)=n$. 

\begin{lemma}
There exists such a vector field $X_n$ for each $n$ in $\NN\setminus\{0\}$.
\end{lemma}
\begin{proof}
Let $\tau$ be a trivialization of the tangent bundle of the three-sphere such that $H_{1,n}$ is a constant map from $S^3$ to $S^2$. Let $T_1$ (resp. $T_2$) be a tubular neighbourhood of $L_1$ (resp. $L_2$), and define $X_n(x)=H_{n,1}(x)$ if $x \in L_1$ and $X_n(x)=-H_{n,1}(x)$ if $x \in L_2$. Let $\star$ be the point of $S^2$ such that $H_{n,1}\equiv\star$ and let $\dagger$ denote a different point of $S^2$. 

Let $\gamma_t,\ t\in [0,1]$ be a smooth path on the two-sphere, such that $\gamma^{-1}(\star)=\{0\}$, $\gamma^{-1}(\dagger)=\{1\}$ and $\gamma_t\neq -\star$ for all $t$.
The solid torus $T_1$ is diffeomorphic to $D^2\times S^1$, where $D^2=\{(r,\theta),0\leq r\leq 1, \theta\in [0,2\pi]\}$. We define $X_n$ on $D^2\times\omega$ by $X_n(r,\theta,\omega)=\gamma(r)$. Therefore, $X_n$ is equal to $H_{1,n}$ on $L_1$, it is equal to $\dagger$ on $\partial T_1$, and is transverse to $H_{1,n}$ in $T_1\setminus L_1$. 

Similarly, we define $X_n$ on $T_2$ such that $X_n\equiv \dagger$ on the boundary of $T_2$. Therefore one can complete $X_n$ outside $T_1\cup T_2$ by $X_n|_{S^3\setminus(T_1\cup T_2)}\equiv \dagger$.

With a small perturbation of $X_n$ along the boundaries of $T_1$ and $T_2$, we can assume that $X_n$ is smooth, always transverse to $H_{n,1}$ except on $L_1$ and $L_2$ where the vector field are positively and negatively tangents, respectively.
\end{proof}

We have therefore a vector field in each homotopy class: $\cI(X_n)=n-1$ for $n>0$.

\section{Homotopy to Morse-Smale flows}
In this section, we will exhibit a non-singular Morse-Smale vector field in each homotopy class of vector fields on the sphere. This result has already been obtained by Yano,~\cite{ya2} after previous work of Wilson,~\cite{wi2}. Yano constructs a non-singular Morse-Smale vector field $X_n$ such that $\cI(X_n)=n$ for each $n$, and $X_n$ has exactly $2n$ periodic orbits. The main interest of our construction is that all the vector fields (except those in the homotopy class of $\cH_+$ and $\cH_-$) have exactly 4 periodic orbits. This shows therefore that the number of those orbits is not relevant (as suspected in~\cite[Remark 5.2]{ya2}). Using our criterion, we see that their linking number is important.

\begin{definition}
A non-singular vector field is a {\em Morse-Smale} vector field if its flow has
\begin{itemize}
\item a finite set of periodic orbits which are all hyperbolic, 
\item the intersections of the invariant manifolds of those orbits are transversal and,
\item its nonwandering set consists entirely of those orbits.
\end{itemize}
\end{definition}

A periodic orbit of a Morse-Smale vector field is either an attractor (index 0), a saddle (index 1) or a repellor (index 2). Following~\cite{mo,as}, Wada, \cite{wa}, has classified all the indexed links that are realisable as the set of periodic orbits of a non-singular Morse-Smale flow on the three-sphere. 

We present here a well-known construction of a Morse-Smale flow in the homotopy class of $\cH_+$, for ours will be similar in other homotopy classes.

\begin{figure}[htb]
\centerline{\includegraphics[scale=.75]{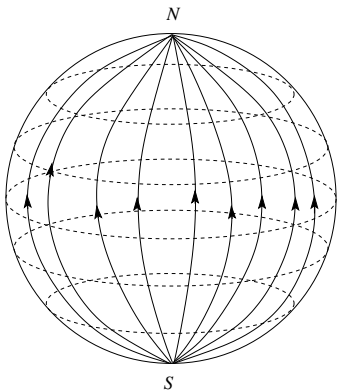}}
\caption{Phase portrait of $X_0$\label{f.sphereNS}}
\end{figure}

Let $X_0$ be the $S-N-$gradient vector field pictured in Figure~\ref{f.sphereNS}.

The Hopf fibration $\cH$ is smooth on $S^3$ and each point of $S^2$ is a regular value of $\cH$. We can lift $X_0$ to a vector field of $S^3$, orthogonal to the orbits of $\cH_+$. The vector field $X_0+\cH_+$ is a Morse-Smale vector field around the $(0,2)-$Hopf link, homotopic to $\cH_+$. The image by $R$ of this vector field gives a Morse-Smale vector field homotopic to $\cH_-$.

We will construct explicitly a Morse-Smale vector field in each homotopy class on the three-sphere. We remark that our vector fields can be obtained using once the fifth or the fourth Wada operation on the $(0,2)-$Hopf link,~\cite{wa}. We will make use of the following remark in our construction.

\begin{remark}
For a given non-singular Morse-Smale vector field, it is always possible to change the orientation of a periodic orbit, staying in the Morse-Smale class.
\end{remark}

The $(p,q)-$Seifert fibration, $\cS_{p,q}$, is smooth on the three-sphere, except along the singular fibres $L_N$ and $L_S$: $\cS_{p,q}$ is smooth on   $\cS_{p,q}^{-1}(S^2\setminus\{N,S\})$ where $N$ and $S$ stand for the North and South poles respectively. Every point of $S^2$, different from $N$ and $S$ is a regular value for $\cS_{p,q}$. This fibration maps the tori $T_{a,b}$ on circles parallel to the equateur of $S^2$.

\begin{figure}[htb]
\centerline{\includegraphics{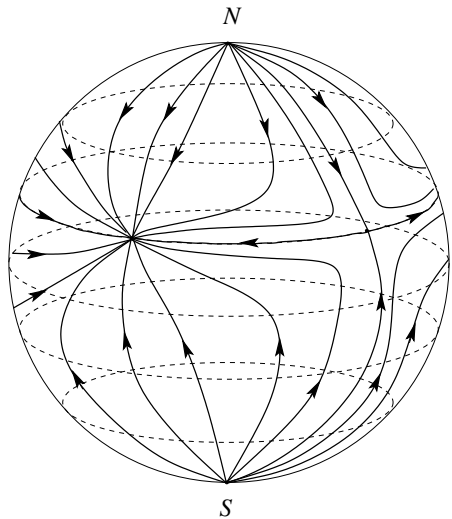}}
\caption{Phase portrait of $X_1$\label{f.sphere}}
\end{figure}

Let $X_1$ be a vector field with two sources (at the poles), one sink and one saddle (that lie on the same parallel), $X_1$ is tangent to this circle (2 heteroclinic orbits) and transverse to every other parallel circles.
The phase portrait of $X_1$ is given on Figure~\ref{f.sphere}.

Let $p$ be a strictly positive integer and let us consider the following vector field $Y_p$ on the three-sphere. For each point in $S^3\setminus\{L_N,L_S\}$, one can lift the vector field $X_1$ on the three-sphere, to a vector field $A$ orthogonal to the fibres of $\cS_{p,1}$. We complete $A$ into a smooth vector field on $S^3$ by $A|_{L_N}=A|_{L_S}=0$. We define $Y_p(x)=H_{p,1}(x)+A(x)$. The vector field $Y_p$ is a Morse-Smale vector field with 4 periodic orbits (corresponding to the singularities of $X_1$), it is positively collinear to $H_{p,1}$ along those orbits, and transverse everywhere else. In fact, it is transverse to every tori $T_{a,b}$ expect the invariant one. On the invariant torus, it is not difficult to make it transverse to $H_{p,1}$ (on the complement of the periodic orbits).

Let $\cM_{2p+1}$ be a Morse-Smale vector field, with the same periodic orbits as $Y_p$ but with the opposite orientation along the orbit of index 0. A slight modification of $Y_p$ leads to $\cM_{2p+1}$, keeping the transversality with $H_{p,1}$ on the complement of the periodic orbits. Similarly, let $\cM_{2p+2}$ be a Morse-Smale vector field, with the same periodic orbits as $\cM_{2p+1}$ but with the opposite orientation along the orbit corresponding to the North pole of the two-sphere.

\begin{figure}[htb]
   \begin{minipage}[c]{.48\linewidth}
      \includegraphics[scale=.85]{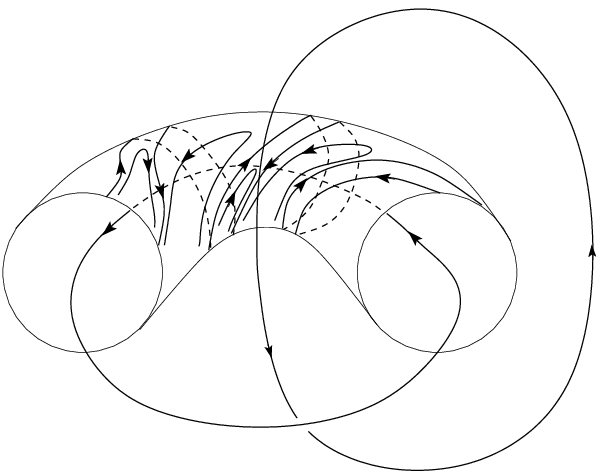}
   \end{minipage} \hfill
   \begin{minipage}[c]{.48\linewidth}
      \includegraphics[scale=.8]{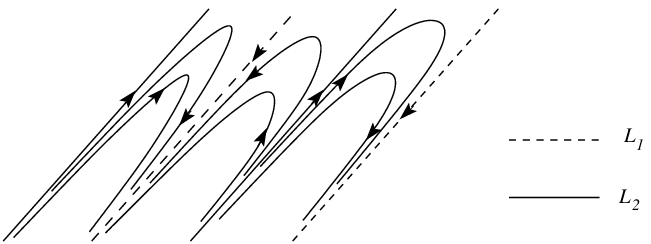}
   \end{minipage}
\caption{The vector field $\cM_{n}$ in $S^3$ and on the invariant torus\label{f.M}}
\end{figure}

As $p$ is strictly positive, we obtain a vector field $\cM_{n}$ for each $n\geq 3$. We define $\cM_2$ as follow (the idea is to do the same with "$H_{0,1}$"): it has 4 periodic orbits, two repellors form a Hopf link, one saddle and one attractor that are linked to one of the repellors, with linking number 1 for the saddle and -1 for the attractor. The saddle and the attractor are on an invariant torus $T_{a,b}$, such that $T_{a,b}$ is the union of the stable manifold of the saddle and the attractor (in particular, the saddle and the attractor are not linked together).

\begin{remark}
The vector field $\cM_2$ is obtained using the fourth Wada operation,~\cite{wa}, and is similar to the vector field constructed by Yano in~\cite{ya2}.

All the other vector fields, $\cM_{n}$ with $n\geq 3$ are obtained using once the fifth Wada operation. 
\end{remark} 

Finally, we define $\cM_1$ to be the Morse-Smale vector field constructed around the $(0,2)-$Hopf link (see the begining of this section).

\begin{lemma}
For $n\geq 1$, the vector field $\cM_{n}$ belongs to the $(n-1)^{th}$ class of homotopy of non-singular vector fields on the three-sphere: for $n\geq 1$, we have $\cI(\cM_{n})=n-1$.
\end{lemma}

\begin{proof}
The lemma is true for $n=1$ by construction.

If $n$ is equal to $2$, $\cM_2$ is tangent to every $H_{p,q}$ with $p$ and $q$ strictly positive, along a Hopf link and two loops on the invariant torus, parallel to the periodic orbits of $\cM_2$. On one of those loops, $H_{p,q}$ and $\cM_2$ are negatively tangent. Then we obtain: $|\mbox{link}(C_+,C_-)|=1$ and $\cI(\cM_2)=1$.

If $n\geq 3$ is an odd integer, $n=2p+1$, the set $C_+$, where $\cM_n$ and $H_{p,1}$ are positively collinear, is the union of the orbits of $\cM_n$ of index 1 and 2, namely $L_N$, $L_S$ and $L_1$. The set $C_-$ is the orbit of index 0: $L_0$. We oriente those orbits with respect to $H_{p,1}$ and we obtain that:
\begin{itemize}
\item $\mbox{link}(L_0,L_1)=\mbox{link}(L_0,L_S)=p$,
\item $\mbox{link}(L_0,L_N)= 1$.
\end{itemize} 
Therefore, with this orientation, $\mbox{link}(C_-,C_+)=2p+1$. By lemma~\ref{l.critere} and lemma~\ref{l.plusmoins}, we have $\cD(\cM_n,\cH_+)=\cD(\cM_n,H_{p,1})=2p+1$.

Similarly, if $n\geq 4$ is even, $n=2p+2$, we have: $C_+$ is the union of $L_S$ and $L_1$, and $C_-$ is the union of $L_N$ and $L_0$. Once again, we oriente those orbits with respect to $H_{p,1}$ and we obtain that:
\begin{itemize}
\item $\mbox{link}(L_S,L_0)=\mbox{link}(L_1,L_0)=p$,
\item $\mbox{link}(L_1,L_N)=\mbox{link}(L_S,L_N)= 1$.
\end{itemize} 
We obtain that $\cD(\cM_n,\cH_+)=2p+2$.

Using the same technique and the picture of $\cM_n$ on the invariant torus, we obtain that $\cD(\cM_{2p+1},\cH_-)=2p$ and $\cD(\cM_{2p},\cH_-)=2p+1$.
Therefore, for each $p>0$, we have $\cI(\cM_{2p+1})=2p=n-1$ and $\cI(\cM_{2p+2})=2p+1=n-1$, then for each $n\geq 1$, $\cI(\cM_n)=n-1$.
\end{proof}

\begin{acknowledgment}
The author would like to thank Robert MacKay, who proposed the project, for many useful conversations, as well as Sebastian van Strien for his constant support.
\end{acknowledgment}

\end{article}
\end{document}